\documentclass{article}

\usepackage[titletoc,title]{appendix}
\usepackage[cp1250]{inputenc}
\usepackage[IL2]{fontenc}
\usepackage{yfonts,fancyhdr}
\usepackage{a4wide}
\usepackage[english]{babel}
\usepackage{euscript}
\usepackage{amstext,amsbsy,amscd,amssymb}
\usepackage{amsmath,enumerate}
\usepackage{amsfonts}
\usepackage{mathrsfs,tensor,xy}
\usepackage{graphics}
\usepackage{microtype}
\include{diagxy}

 \input xypic
 \input epsf
 \input xy
 \xyoption{all}

\let\rarr=\rightarrow

\let\veps=\varepsilon
\let\mcal=\mathcal
\let\mfrak=\mathfrak
\let\eus=\EuScript

\def\N{\mathbb{N}}

\def\R{\mathbb{R}}
\def\C{\mathbb{C}}

\def\End{\mathop {\rm End} \nolimits}
\def\Hom{\mathop {\rm Hom} \nolimits}

\def\Ind{\mathop {\rm Ind} \nolimits}

\def\ad{\mathop {\rm ad} \nolimits}

\def\Diff{\mathop {\rm Diff} \nolimits}

\def\GL{\mathop {\rm GL} \nolimits}

\def\SO{\mathop {\rm SO} \nolimits}
\def\Spin{\mathop {\rm Spin} \nolimits}

\def\Sol{\mathop {\rm Sol} \nolimits}

\def\tr{\mathop {\rm tr} \nolimits}

\newcommand{\mP}{\mathbb P}

\newcommand{\mS}{\mathbb S}
\newcommand{\mV}{\mathbb V}

\newcommand{\mC}{\mathbb C}

\newcommand{\mN}{\mathbb N}

\newcommand{\lC}{\cal C}

\newcommand{\mR}{\mathbb R}
\newcommand{\mW}{\mathbb W}

\newcommand{\gog}{\mathfrak{g}}

\newcommand{\gop}{\mathfrak{p}}

\newsymbol\squares 1003

\long\def\proof #1{\noindent \emph{Proof.}\ #1 \hfill $\squares$
\medskip}

\newcounter{num}[section]
\numberwithin{equation}{section}
\numberwithin{num}{section}

\long\def\definition #1 {\refstepcounter{num} \noindent {\bf
Definition \thenum.} #1

\medskip}

\long\def\theorem #1{\refstepcounter{num} \noindent {\bf Theorem
\thenum.} #1

\medskip}

\long\def\lemma #1{\refstepcounter{num}  \noindent {\bf Lemma
\thenum.} #1

\medskip}

\newcommand*\riso{%
  \xrightarrow[]{\raisebox{-0.25em}{\smash{\ensuremath{\sim}}}}%
}

\makeatletter
\newcommand*\if@single[3]{%
  \setbox0\hbox{${\mathaccent"0362{#1}}^H$}%
  \setbox2\hbox{${\mathaccent"0362{\kern0pt#1}}^H$}%
  \ifdim\ht0=\ht2 #3\else #2\fi
  }
\newcommand*\rel@kern[1]{\kern#1\dimexpr\macc@kerna}
\newcommand*\widebar[1]{\@ifnextchar^{{\wide@bar{#1}{0}}}{\wide@bar{#1}{1}}}
\newcommand*\wide@bar[2]{\if@single{#1}{\wide@bar@{#1}{#2}{1}}{\wide@bar@{#1}{#2}{2}}}
\newcommand*\wide@bar@[3]{%
  \begingroup
  \def\mathaccent##1##2{%
    \if#32 \let\macc@nucleus\first@char \fi
    \setbox\z@\hbox{$\macc@style{\macc@nucleus}_{}$}%
    \setbox\tw@\hbox{$\macc@style{\macc@nucleus}{}_{}$}%
    \dimen@\wd\tw@
    \advance\dimen@-\wd\z@
    \divide\dimen@ 3
    \@tempdima\wd\tw@
    \advance\@tempdima-\scriptspace
    \divide\@tempdima 10
    \advance\dimen@-\@tempdima
    \ifdim\dimen@>\z@ \dimen@0pt\fi
    \rel@kern{0.6}\kern-\dimen@
    \if#31
      \overline{\rel@kern{-0.6}\kern\dimen@\macc@nucleus\rel@kern{0.4}\kern\dimen@}%
      \advance\dimen@0.4\dimexpr\macc@kerna
      \let\final@kern#2%
      \ifdim\dimen@<\z@ \let\final@kern1\fi
      \if\final@kern1 \kern-\dimen@\fi
    \else
      \overline{\rel@kern{-0.6}\kern\dimen@#1}%
    \fi
  }%
  \macc@depth\@ne
  \let\math@bgroup\@empty \let\math@egroup\macc@set@skewchar
  \mathsurround\z@ \frozen@everymath{\mathgroup\macc@group\relax}%
  \macc@set@skewchar\relax
  \let\mathaccentV\macc@nested@a
  \if#31
    \macc@nested@a\relax111{#1}%
  \else
    \def\gobble@till@marker##1\endmarker{}%
    \futurelet\first@char\gobble@till@marker#1\endmarker
    \ifcat\noexpand\first@char A\else
      \def\first@char{}%
    \fi
    \macc@nested@a\relax111{\first@char}%
  \fi
  \endgroup
}
\makeatother

\newcommand\rsmraise[1]{%
  \ifx#1\displaystyle .8\else
    \ifx#1\textstyle .8\else
      \ifx#1\scriptstyle .6\else
        .45%
      \fi
    \fi
  \fi}

\title{Equivariant differential operators on spinors in conformal geometry}

\author{Libor Křižka, Petr Somberg}

\AtEndDocument{\bigskip{\footnotesize%
  (L.\,Křižka) \textsc{Charles University, Faculty of Mathematics and Physics,
	Mathematical Institute, \\ Sokolovská 83, 180\,00 Praha 8, Czech Republic} \par
  \textit{E-mail address}: \texttt{krizka@karlin.mff.cuni.cz} \par
  \addvspace{\medskipamount}
  (P.\,Somberg) \textsc{Charles University, Faculty of Mathematics and Physics,
	Mathematical Institute, \\ Sokolovská 83, 180\,00 Praha 8, Czech Republic}  \par
  \textit{E-mail address}: \texttt{somberg@karlin.mff.cuni.cz} \par
}}

\begin{document}
\date{}
\maketitle

\begin{abstract}
We present a novel approach to the classification of conformally equivariant
differential operators on spinors in the case of homogeneous conformal geometry.
It is based on the classification of solutions for a vector-valued system
of partial differential equations, associated to $\mathcal{D}$-modules for the homogeneous conformal
structure and controlled by the spin Howe duality for the orthogonal Lie algebras.


\end{abstract}

\thispagestyle{empty}

\tableofcontents


\section*{Introduction}
\addcontentsline{toc}{section}{Introduction}

The concept of conformal symmetry canonically extends (pseudo)Riemannian structure, and
so it naturally arises in various functional, analytical and
geometrical problems on manifolds equipped with a conformal class of metrics.

There are two types of conformally equivariant differential operators acting on
spinor fields. Namely, in the case of positive definite metrics we have a class of elliptic
operators called the conformal powers of the Dirac operator and a class of
overdetermined twistor operators. The basic first order representatives in these sets, the Dirac and
the twistor operators, are ubiquitous in analysis, geometry and representation
theory, see e.g.\ \cite{Friedrich2000}, \cite{Baum-Friedrich-Grunewald-Kath1991},
\cite{Fischmann2014}, and the references therein.

In the present short article we exploit the general framework of \cite{KOSS},
\cite{Krizka-Somberg2015}, \cite{Kostant1975}, and give a short, self-content
and signature of the underlying metric independent classification statement,
which can be easily modified and adopted
to many analogous problems (e.g.\ the classification for higher spinor representations).
This classification result concerns the conformally equivariant differential operators
acting on spinor fields. Though well-known to specialists, it is difficult
to find an explicit statement in the existing literature. The first such classification
was presented in the paper \cite{Boe-Collingwood1985}. For a general classification
scheme we refer to \cite{Slovak1992}, but being expressed in terms of the representation
theoretical data it lacks the explicit description of the origin and presentation
for such operators. The translation principle was used to obtain the classification
of curved extensions in \cite{Eastwood-Slovak1997}.

Our approach is based on the techniques of $\mcal{D}$-modules, emerging
in the reformulation of the former problem for the existence and construction
of conformally equivariant differential operators to a question on the
existence and construction of the space of homomorphisms between certain
algebraic objects called conformal generalized Verma modules. The space
of homomorphisms is detected by special
elements called singular vectors, which are the solution spaces of a system of
PDEs constructed out of the action of certain elements in the conformal Lie
algebra. This observation in tandem with the spin Howe duality for the
orthogonal Lie algebra (the simple part of the Levi subalgebra of the conformal
Lie algebra) allows to complete the classification task.

The content of our article goes as follows.
Section \ref{sec:Verma modules} contains a general introduction into
the techniques used in the article, based on the approach of
algebraic analysis on $\mcal{D}$-modules and leading to the class of conformal
generalized Verma modules needed throughout the letter.
In Section \ref{sec:conformal structure} we overview several geometrical
and representation theoretical aspects of the homogeneous conformal
structure with emphasis on the spinor representation, and then in
Section \ref{sec:singular vectors} formulate
and prove the classification result for the singular vectors in conformal
generalized Verma modules induced from the spinor representation twisted
by characters. In Section \ref{sec:invariant operators} we dualize the results achieved in Section
\ref{sec:singular vectors} and obtain the classification of conformal equivariant differential
operators on spinors. For the reader's convenience, we summarize the
spin Howe duality for the orthogonal Lie algebra
in terms of the Fischer decomposition.

Due to a uniform description of the conformal symmetry we restrict to
the underlying dimension $n \geq 3$. Throughout
the article, $\N$ denotes the natural numbers and $\N_0$
denotes the natural numbers including zero.


\section{Equivariant differential operators and algebraic analysis on generalized Verma modules}
\label{sec:Verma modules}

Let us consider the pair $(G,P)$, consisting of a connected real reductive Lie group $G$ and
its parabolic subgroup $P$. In the Levi decomposition $P=LU$,
$L$ denotes the Levi subgroup and $U$ the unipotent radical of $P$.
We write $\mfrak{g}(\R),\,\mfrak{p}(\R),\,\mfrak{l}(\R),\,\mfrak{u}(\mR)$ for the
real Lie algebras and $\mfrak{g},\,\mfrak{p},\,\mfrak{l},\,\mfrak{u}$ for the complexified Lie
algebras of $G,\,P,\,L,\,U$, respectively. The symbol $U$ applied to a Lie
algebra denotes its universal enveloping algebra.

It is well-known that the $G$-equivariant differential operators acting on
principal series representations for $G$ can be recognized in the study
of homomorphisms between generalized Verma modules for the Lie algebra
$\mfrak{g}$. The latter homomorphisms are determined by the image of the highest
weight vectors and are referred to as the singular vectors, characterized as the
vectors in the generalized Verma module annihilated by the positive nilradical $\mfrak{u}$.

We shall rely on the following approach to find precise positions of singular
vectors in a given representation space, cf.\ \cite{KOSS}, \cite{Krizka-Somberg2015}
for a detailed exposition. First of all, in the present article $\mV$ denotes a complex simple
finite dimensional $L$-module, extended to $P$-module by $U$ acting trivially. We denote
by $\mV^*$ the dual $P$-module to $\mV$. Let us assume that
$\lambda \in \Hom_P(\mfrak{p},\C)$ defines a group character
$e^\lambda \colon P \rarr \GL(1,\C)$ of $P$, and define
$\rho \in \Hom_P(\mfrak{p},\C)$ by
\begin{align}
  \rho(X) = {\textstyle {1 \over 2}} \tr_\mfrak{u}\ad(X) \label{eq:rho vector}
\end{align}
for $X \in \mfrak{p}$. Then we define a twisted $P$-module
$\mathbb{V}_{\lambda+\rho}$ with a twist $\lambda +\rho \in \Hom_P(\mfrak{p},\C)$,
where $p \in P$ acts as $e^{\lambda+\rho}(p)p\cdot v$ instead of
$p\cdot v$ for all $v \in \mathbb{V}_{\lambda+\rho} \simeq \mathbb{V}$ (the isomorphism of vector spaces).
In the present article, $\mV=\mathbb{S}$ is the
finite-dimensional (semi)simple $L \simeq (\GL(1,\R)_+\! \times \Spin(p,q,\R))$-module
($p+q \geq 3$) given by the spinor representation twisted by character.

For a chosen principal series representation of $G$ on the vector space $\Ind_P^G(\mV_{\lambda+\rho})$ of smooth sections of the homogeneous vector bundle
$G \times_P \mV_{\lambda+\rho} \rarr G/P$ associated to a $P$-module $\mV_{\lambda+\rho}$, we compute the infinitesimal action
\begin{align}
\pi_\lambda \colon \mfrak{g} \rarr \mcal{D}(U_e) \otimes_\C \End \mathbb{V}_{\lambda+\rho}.
\end{align}
Here $\mcal{D}(U_e)$ denotes the $\C$-algebra of smooth complex linear differential operators on $U_e=\widebar{U}P \subset G/P$ ($\widebar{U}$ is the Lie group whose Lie algebra is the opposite nilradical $\widebar{\mathfrak{u}}(\R)$ to $\mathfrak{u}(\R)$), on the vector space $\mcal{C}^\infty(U_e) \otimes_\C\! \mathbb{V}_{\lambda+\rho}$ of $\mathbb{V}_{\lambda+\rho}$-valued smooth functions on $U_e$ in the non-compact picture of the induced representation.

Since the vector space $\mcal{D}'_o(U_e)\otimes_\C\! \mV_{\lambda+\rho}$ of $\mathbb{V}_{\lambda+\rho}$-valued distributions on $U_e$ supported in the unit coset
$o=eP\in G/P$ is $\mcal{D}(U_e) \otimes_\C \End \mathbb{V}_{\lambda+\rho}$-module, we obtain the infinitesimal action of $\pi_\lambda(X)$ for $X \in \mfrak{g}$ on $\mcal{D}'_o(U_e)\otimes_\C\! \mV_{\lambda+\rho}$.
The exponential map allows to identify $U_e$ with the nilpotent Lie algebra
$\widebar{\mfrak{u}}(\R)$. Denoting by $\eus{A}^\mfrak{g}_{\widebar{\mfrak{u}}}$
the Weyl algebra of the complex vector space $\widebar{\mfrak{u}}$, the vector
space $\mcal{D}'_o(U_e)$ can be conveniently analyzed by identifying it as an
$\eus{A}^\mfrak{g}_{\widebar{\mfrak{u}}}$-module with the quotient of
$\eus{A}^\mfrak{g}_{\widebar{\mfrak{u}}}$ by the left ideal $I_e$ generated
by all polynomials on $\widebar{\mathfrak{u}}$ vanishing at the origin.
Moreover, there is a $U(\gog)$-module isomorphism
\begin{align}\label{vermadistr}
\Phi_\lambda \colon M^\gog_\gop(\mV_{\lambda-\rho})\equiv U(\gog)\otimes_{U(\gop)}\!\mV_{\lambda-\rho} \rarr  \mcal{D}'_o(U_e)\otimes_\C\! \mV_{\lambda+\rho} \simeq \eus{A}^\mfrak{g}_{\widebar{\mfrak{u}}}/I_e \otimes_\C\! \mathbb{V}_{\lambda+\rho}.
\end{align}
Let $(x_1,x_2,\dots,x_n)$ be the linear coordinate functions on $\widebar{\mfrak{u}}$ and let $(y_1,y_2,\dots,y_n)$ be the dual linear coordinate functions on $\widebar{\mfrak{u}}^*$. Then the algebraic Fourier transform
\begin{align}
  \mcal{F} \colon \eus{A}^\mfrak{g}_{\widebar{\mfrak{u}}} \rarr \eus{A}^\mfrak{g}_{\widebar{\mfrak{u}}^*}
\end{align}
is given by
\begin{align}
  \mcal{F}(x_i) = -\partial_{y_i}, \qquad \mcal{F}(\partial_{x_i}) = y_i
\end{align}
for $i=1,2,\dots,n$, and gives a vector space isomorphism
\begin{align}\label{eqn:FT}
\tau \colon  \eus{A}^\mfrak{g}_{\widebar{\mfrak{u}}}/I_e  \riso  \eus{A}^\mfrak{g}_{\widebar{\mfrak{u}}^*}/\mcal{F}(I_e) \simeq \C[\widebar{\mfrak{u}}^*]
\end{align}
defined by
\begin{align}
  Q\ {\rm mod}\ I_e \mapsto \mcal{F}(Q)\ {\rm mod}\ \mcal{F}(I_e)
\end{align}
for $Q \in \eus{A}^\mfrak{g}_{\widebar{\mfrak{u}}}$. The composition of the previous mappings \eqref{vermadistr} and \eqref{eqn:FT} gives a vector space isomorphism
\begin{align}\label{eqn:phiPM}
 \tau \circ \Phi_\lambda \colon U(\mfrak{g}) \otimes_{U(\mfrak{p})} \! \mathbb{V}_{\lambda-\rho} \riso \mcal{D}'_o(U_e) \otimes_\C \!\mV_{\lambda+\rho} \riso
\C[\widebar{\mathfrak{u}}^*]\otimes_\C\! \mV_{\lambda-\rho},
\end{align}
thereby inducing the action $\hat{\pi}_\lambda \colon \mfrak{g} \rarr \eus{A}^\mfrak{g}_{\widebar{\mfrak{u}}^*}\! \otimes_\C \End \mathbb{V}_{\lambda-\rho}$ of $\gog$ on $\C[\widebar{\mathfrak{u}}^*]\otimes_\C\! \mV_{\lambda-\rho}$.
The polynomial algebra on $\widebar{\mfrak{u}}^*$ is isomorphic to the universal enveloping algebra
$U(\widebar{\mathfrak{u}})$.
\medskip

\definition{Let $\mV$ be a complex (semi)simple finite-dimensional $L$-module,
extended to a $P$-module by $U$ acting trivially. We define the $L$-module
\begin{align}
M_\mfrak{p}^\mfrak{g}(\mV)^\mfrak{u}= \{v\in M^\mfrak{g}_\mfrak{p}(\mV);\, Xv=0\
\text{for all}\ X\in \mathfrak{u}\}
\end{align}
and denote it the vector space of singular vectors.}

The vector space of singular vectors is for any finite-dimensional
complex (semi)simple $P$-module $\mV$ a finite-dimensional
completely reducible $L$-module. We denote by $\mW$ one of its simple
$L$-submodules, and this yields
$U(\gog)$-homomorphism from $M_\gop^\gog(\mW)$ to $M_\gop^\gog(\mV)$ such that
\begin{align}
\Hom_{(\gog,P)}(M_\gop^\gog(\mW), M_\gop^\gog(\mV))
\simeq
\Hom_L(\mW, M_\gop^\gog(\mV)^\mfrak{u}).
\end{align}
We introduce the $L$-module
\begin{align} \label{eqn:sol2l}
\Sol(\mathfrak{g},\mathfrak{p};\C[\widebar{\mathfrak{u}}^*]\otimes_\C\! \mV_{\lambda-\rho})^\mcal{F} =\{f \in \C[\widebar{\mathfrak{u}}^*]\otimes_\C\! \mV_{\lambda-\rho};\, \hat{\pi}_\lambda(X) f = 0\ \text{for all}\ X\in\mathfrak{u}\},
\end{align}
and by \eqref{eqn:phiPM}, there is an $L$-equivariant isomorphism
\begin{equation}
\label{eqn:phi}
\tau \circ \Phi_\lambda \colon M_\mathfrak{p}^\mathfrak{g}(\mV_{\lambda-\rho})^\mathfrak{u} \riso \Sol(\mathfrak{g},\mathfrak{p}; \C[\widebar{\mathfrak{u}}^*]\otimes_\C\! \mV_{\lambda-\rho})^\mcal{F}.
\end{equation}
The action of $\hat{\pi}_\lambda(X)$ on $\C[\widebar{\mathfrak{u}}^*]\otimes_\C\! \mV_{\lambda-\rho}$
produces a system of partial differential equations for the elements in
$\Sol(\mathfrak{g},\mathfrak{p}; \C[\widebar{\mathfrak{u}}^*]\otimes_\C\! \mV_{\lambda-\rho})^\mcal{F}$, which makes possible to describe completely the structure of its solution space in particular cases of interest. Namely, the algebraic Fourier transform on $\mcal{D}$-modules converts the algebraic problem of finding singular vectors in generalized Verma
modules into an analytic problem of solving the systems of partial differential equations.

The formulation above has the following classical dual statement,
cf.\ \cite{Collingwood-Shelton1990}, \cite{Kostant1975}, which
explains the relationship between the geometrical problem of finding
$G$-equivariant differential operators between induced representations
and the algebraic problem of finding homomorphisms between generalized
Verma modules. Let $\mV$ and $\mW$ be two simple finite-dimensional $P$-modules.
Then the vector space of $G$-equivariant differential operators
$\Hom_{\Diff(G)}(\Ind_P^G(\mV),\Ind_P^G(\mW))$
is isomorphic to the vector space of $(\gog,P)$-homomorphisms
$\Hom_{(\gog,P)}(M^\gog_\gop(\mW^*),M^\gog_\gop(\mV^*))$.

\section{Geometrical and representation theoretical aspects of\\ conformal structure}
\label{sec:conformal structure}

In the present section we describe the rudiments of the geometry
of the homogeneous (flat) conformal structure, with emphasis on
representation theoretical aspects of the spinor representation.
The dimension $n$ of the underlying space is supposed to be at
least three due to a uniform behavior of the conformal symmetry for
$n\geq 3$. We shall treat the case of general signature $(p,q)$,
$p+q=n$, so the conformal Lie group is the real form $\SO_o(p+1,q+1,\mR)$.

\subsection{The geometrical model of homogeneous conformal structure}

The generalized flag manifold describing the homogeneous (flat) model of
real conformal structure of signature $(p,q)$, $p+q=n$, is the homogeneous space
$G/P\simeq S^{p,q}$, where $S^{p,q}$ is the quadric in $\R^{p+1,q+1}$ and the
connected Lie group $G=\SO_o(p+1,q+1,\mR)$ is
the group of automorphisms of the vector space $\R^{p+1,q+1}$ preserving
the inner product $\langle\cdot\,,\cdot\rangle$
of signature $(p+1,q+1)$ on $\R^{p+1,q+1}$ corresponding to the matrix
\begin{align}
  J_{p+1,q+1} =
  \begin{pmatrix}
    0 & 0 & 1 \\
    0 & I_{p,q} & 0 \\
    1 & 0 & 0
  \end{pmatrix}\!,
  \qquad
  I_{p,q}=
  \begin{pmatrix}
    I_p & 0 \\
    0 & -I_q
  \end{pmatrix}\!,
\end{align}
where $I_n$ means the identity matrix $n \times n$, and $P\subset G$ is the conformal parabolic subgroup.

Let $(x_0,x_1,\dots,x_n,x_\infty)$ be the canonical linear coordinate functions on $\R^{p+1,q+1}$,
and denote $x_A=(x_0,x_1,\dots,x_n,x_\infty)=(x_0,x_a,x_\infty)$
with $x_a=(x_1,x_2,\dots,x_n)$. The null cone
of $(\R^{p+1,q+1}, \langle\cdot\,,\cdot\rangle)$ is defined as
$\mcal{N}_{p,q}=\{x_A\in \R^{p+1,q+1};\, \langle x_A,x_A\rangle=0\}$,
and its projectivization $\mP(\mcal{N}_{p,q})$ is the conformal sphere $(S^{p,q},[g_0])$. The conformal class $[g_0]$
of the round metric $g_0$ on $S^{p,q}$ is induced by the map
$x_a\rarr  \langle x_A\, ,x_A\rangle$, where $x_a\in TS^{p,q}$. The tangent space
$TS^{p,q}$ is identified with the space of classes of vector fields on $\R^{p+1,q+1}$
of homogeneity one, which are orthogonal to $x_A$ with
respect to $\langle\cdot\,,\cdot\rangle$. A $\lambda$-density on $S^{p,q}$ is identified
with a function of homogeneity $\lambda$ on the null cone
$\mcal{N}_{p,q}$. As for the spinor bundle $\mcal{S}$ on $S^{p,q}$,
it can be identified with the associated vector bundle
whose fiber at the point on $S^{p,q}$ corresponding to the line
$[v]\in \mcal{N}_{p,q}$ is given by the quotient space
$\tilde{\mcal{S}}/({x_A}\cdot\tilde{\mcal{S}})$.
Here we used the notation $\tilde{\mcal{S}}$ for the trivial spinor bundle on $\R^{p+1,q+1}$, and the vector
$x_A$ is regarded as an element of $\End(\tilde{\mcal{S}})\simeq \mcal{C}\ell(\mR^{p+1,q+1}, J_{p+1,q+1})$
acting by the Clifford multiplication on $\tilde{\mcal{S}}$.

The Lie group $G$ acts transitively on the space of lines in
$\mcal{N}_{p,q}$ by $(g,[v]) \mapsto [g.v]$ for $0\not=v\in\mR^{p+1,q+1}$,
$g\in\SO_o(p+1,q+1,\mR)$. The stabilizer of the null ray $[(1,0,\dots,0)]$
is the real conformal parabolic subgroup
$P \simeq (\GL(1,\R)_+ \times \SO(p,q,\R))\ltimes \R^{p,q}$,
\begin{align}
 P=\left\{\!
  \begin{pmatrix}
    \chi(A)a & -av^{\rm t}I_{p,q}A & -{a\over 2} v^{\rm t}I_{p,q}v \\
    0 & A & v \\
	0 & 0 & \chi(A)a^{-1}
  \end{pmatrix}\!;
  \begin{gathered}
  a\in\R_+,\, v\in\R^{p,q},\, A \in M_{n,n}(\R)\\ \det A =1,\, A^{\rm t}I_{p,q}A=I_{p,q}
  \end{gathered}\right\}\!,
\end{align}
where the superscript in $v^{\rm t}$ denotes the transpose of
$v$ and $\chi \colon \SO(p,q,\R) \rarr \{-1,1\}$ is the
multiplicative character defined by $\chi(A)=1$ for $A\in \SO_o(p,q,\R)$
and $\chi(A)=-1$ if $A$ does not belong to the component of identity of $\SO(p,q,\R)$.
In the Dynkin diagrammatic notation, the Lie algebra of $P$ is given
by omitting the first simple root in the $B, D$-series of simple real Lie
algebras.

Because of our interest in the half-integral representations,
we need the double cover of the connected component of the Lie group $G=\SO_o(p+1,q+1,\R)$
called the spin group $\smash{\widetilde{G}}=\Spin_o(p+1,q+1,\R)$.
There is a double cover homomorphism $\psi \colon \smash{\widetilde{G}} \rarr G$ of Lie groups, and
if we define the parabolic subgroup $\smash{\widetilde{P}}$ of
$\smash{\widetilde{G}}$ by $\smash{\widetilde{P}}=\psi^{-1}(P)$, then the simple part
of the Levi subgroup of $\smash{\widetilde{P}}$ is of the form
$\smash{\widetilde{M}}\simeq \Spin(p,q,\R)$
and the mapping $\psi$ induces an isomorphism of generalized flag
manifolds $\smash{\widetilde{G}/\widetilde{P}}\simeq G/P$.

For more detailed introduction into conformal geometry we refer e.g.\ to
\cite{Curry-Gover2014} and the references therein.

\subsection{Representation theory of conformal geometry}
\label{section3.1}

Let us consider the connected complex simple Lie group $G_\C=\SO(n+2,\C)$, $n\geq 3$, defined by
\begin{align}
  \SO(n+2,\C)=\{X\in \GL(n+2,\C);\, X^{\rm t}J_{p+1,q+1}X=J_{p+1,q+1}\},
\end{align}
and its Lie algebra $\mfrak{g}=\mfrak{so}(n+2,\C)$ given by
\begin{align}
\begin{aligned}
  \mfrak{so}(n+2,\C)&=\{X\in M_{n+2,n+2}(\C);\, X^{\rm t}J_{p+1,q+1}+J_{p+1,q+1}X=0\} \\
  &=\left\{\!
  \begin{pmatrix}
    a & v^{\rm t} & 0 \\
    u & A & -I_{p,q}v \\
    0 & -u^{\rm t}I_{p,q} & -a
  \end{pmatrix}\!;
  \begin{gathered}
  a\in\C,\, u,v \in \C^m,\, A \in M_{n,n}(\C),\\ A^{\rm t}I_{p,q}+I_{p,q}A=0
  \end{gathered}\right\}.
\end{aligned}
\end{align}
The standard parabolic subgroup $P_\C$ of $G_\C$ is defined by
\begin{align}
 P=\left\{\!
  \begin{pmatrix}
    a & -av^{\rm t}I_{p,q}A & -{a\over 2} v^{\rm t}I_{p,q}v \\
    0 & A & v \\
	0 & 0 & a^{-1}
  \end{pmatrix}\!;
  \begin{gathered}
  a\in\C,\, v\in\C^n,\, A \in M_{n,n}(\C),\\ a\neq0,\, \det A=1,\, A^{\rm t}I_{p,q}A=I_{p,q}
  \end{gathered}\right\}
\end{align}
and its Lie algebra $\mfrak{p}$ is given by
\begin{align}
  \mfrak{p}=\left\{\!
  \begin{pmatrix}
    a & v^{\rm t} & 0 \\
    0 & A & -I_{p,q}v \\
    0 & 0 & -a
  \end{pmatrix}\!;
  \begin{gathered}
  a\in\C,\, u,v \in \C^m,\, A \in M_{n,n}(\C),\\ A^{\rm t}I_{p,q}+I_{p,q}A=0
  \end{gathered}\right\}.
\end{align}
Let us denote by $\mfrak{u}$ the nilradical of the parabolic subalgebra $\mfrak{p}$ and by $\widebar{\mfrak{u}}$ the opposite niradical. Then we have a triangular decomposition $\mfrak{g} = \widebar{\mfrak{u}} \oplus \mfrak{l} \oplus \mfrak{u}$, where $\mfrak{l}$ is the Levi subalgebra of $\mfrak{p}$.

We choose a basis $(f_1,f_2,\dots,f_n)$ of the commutative opposite nilradical
$\widebar{\mfrak{u}}$ by
\begin{align}
  f_i =
  \begin{pmatrix}
    0 & 0 & 0 \\
    1_i & 0 & 0 \\
    0 & -\veps_i 1_i^{\rm t} & 0
  \end{pmatrix}\!,
\end{align}
where $\veps_i=1$ for $i=1,2,\dots,p$ and $\veps_i=-1$ for $i=p+1,p+2,\dots,n$, and a
basis $(g_1,g_2,\dots,g_n)$ of the commutative nilradical $\mfrak{u}$ by
\begin{align}
  g_i =
  \begin{pmatrix}
    0 & 1_i^{\rm t} & 0 \\
    0 & 0 & -\veps_i1_i \\
    0 & 0 & 0
  \end{pmatrix}\!.
\end{align}
The Levi subalgebra $\mfrak{l}$ of $\mfrak{p}$ is the linear span of
\begin{align}
  h =
  \begin{pmatrix}
    1 & 0 & 0 \\
    0 & 0 & 0 \\
    0 & 0 & -1
  \end{pmatrix}\!, \qquad
  h_A =
  \begin{pmatrix}
    0 & 0 & 0 \\
    0 & A & 0 \\
    0 & 0 & 0
  \end{pmatrix}\!,
\end{align}
where $A \in M_{n,n,}(\C)$ satisfies $A^{\rm t}I_{p,q}+I_{p,q}A=0$. Moreover, the element $h$ is
a basis of the center $\mfrak{z}(\mfrak{l})$ of $\mfrak{l}$.

The real connected simple Lie group $G$ and its real parabolic subgroup $P$ are defined as the identity
components of $G_\C \cap \GL(n+2,\R)$ and $P_\C \cap \GL(n+2,\R)$, respectively, and their real Lie
algebras are $\mfrak{g}(\R)$ and $\mfrak{p}(\R)$, respectively.

Any character $\lambda \in \Hom_P(\mfrak{p},\C)$ is given by
\begin{align}
  \lambda = \alpha \widetilde{\omega}
\end{align}
for some $\alpha \in \C$, where $\widetilde{\omega}\in \Hom_P(\mfrak{p},\C)$
is defined by $\widetilde{\omega}(h)=1$, $\widetilde{\omega}(h_A)=0$ and then
trivially extended to
$\mfrak{l} \oplus \mfrak{u}$. The vector $\rho \in \Hom_P(\mfrak{p},\C)$ defined by \eqref{eq:rho vector} is
\begin{align}
  \rho={\textstyle {n \over 2}} \widetilde{\omega}.
\end{align}
By abuse of notation, we use the simplified notation $\lambda \in \Hom_P(\mfrak{p},\C)$
for the character $\lambda \widetilde{\omega} \in \Hom_P(\mfrak{p},\C)$, $\lambda \in \C$.

\subsection{Description of the representation}

Here we describe the representations of $\mathfrak{g}$
on the space of sections of vector bundles on $\smash{\widetilde{G}/\widetilde{P}} \simeq G/P$ associated to the
(semi)simple spinor representation $\mathbb{S}_{\lambda+\rho}$ of $\smash{\widetilde{P}}$
twisted by characters $\lambda + \rho \in \Hom_P(\mfrak{p},\C)$.

The induced representations in question are described in
the non-compact picture, restricting sections on $G/P$ to the
open Schubert cell $U_e$ isomorphic by the exponential
map to the opposite nilradical $\widebar{\mfrak{u}}(\R)$.
Let us denote by $(\hat{x}_1,\hat{x}_2,\dots,\hat{x}_n)$ the linear coordinate functions on
$\widebar{\mfrak{u}}$ with respect to the basis $(f_1,f_2,\dots,f_n)$
of the opposite nilradical $\widebar{\mfrak{u}}$, and by
$(x_1,x_2,\dots,x_n)$ the dual linear coordinate functions on $\widebar{\mfrak{u}}^*$. Then the Weyl algebra
$\eus{A}^\mfrak{g}_{\widebar{\mfrak{u}}}$ is generated by
\begin{align}
\{\hat{x}_1,\dots,\hat{x}_n,\partial_{\hat{x}_1},\dots,\partial_{\hat{x}_n}\}
\label{eq:Weyl algebra generators}
\end{align}
and the Weyl algebra
$\eus{A}^\mfrak{g}_{\widebar{\mfrak{u}}^*}$ is generated by
\begin{align}
\{x_1,\dots,x_n,\partial_{x_1},\dots,\partial_{x_n}\}.
\label{eq:dual Weyl algebra generators}
\end{align}
The local coordinate chart $u_e \colon x\in U_e \mapsto u_e(x)\in \widebar{\mfrak{u}}(\R) \subset \widebar{\mfrak{u}}$ for the open subset $U_e\subset G/P$, in coordinates with respect to the basis
$(f_1,f_2,\dots,f_n)$ of $\widebar{\mfrak{u}}$, is given by
\begin{align}
  u_e(x)=\sum_{i=1}^n u^i(x)f_i \label{eq:coordinate function}
\end{align}
for all $x \in U_e$.

Let $(\sigma,\mathbb{V})$, $\sigma \colon \mfrak{p} \rarr \mfrak{gl}(\mathbb{V})$, be a $\mfrak{p}$-module.
Then a twisted $\mfrak{p}$-module $(\sigma_\lambda, \mathbb{V}_\lambda)$,
$\sigma_\lambda \colon \mfrak{p} \rarr \mfrak{gl}(\mathbb{V}_\lambda)$, with a twist
$\lambda \in \Hom_P(\mfrak{p},\C)$, is defined as
\begin{align}
  \sigma_\lambda(X)v=\sigma(X)v+\lambda(X)v
\end{align}
for all $X \in \mfrak{p}$ and $v \in \mathbb{V}_\lambda \simeq \mathbb{V}$ (as vector spaces).

Let us introduce the notation
\begin{align}
  E_x = {\textstyle \sum\limits_{j=1}^n} x_j\partial_{x_j} \qquad \text{and} \qquad
  E_{\hat{x}} = {\textstyle \sum\limits_{j=1}^n} \hat{x}_j\partial_{\hat{x}_j}.
\end{align}
for the Euler homogeneity operators.
\medskip

\theorem{\label{reproper}
Let $\lambda \in \Hom_P(\mfrak{p},\C)$ and let
$(\sigma,\mathbb{V})$, $\sigma \colon \mfrak{p} \rarr \mfrak{gl}(\mathbb{V})$,
be a $\mfrak{p}$-module. Then the embedding of $\mfrak{g}$ into
$\eus{A}^\mfrak{g}_{\widebar{\mfrak{u}}} \otimes_\C \End \mathbb{V}_{\lambda+\rho}$ and
$\eus{A}^\mfrak{g}_{\widebar{\mfrak{u}}^*}\! \otimes_\C \End \mathbb{V}_{\lambda-\rho}$
is given by
\begin{enumerate}
\item[1)]
\begin{align}
  \pi_\lambda(f_i)&=-\partial_{\hat{x}_i},
\end{align}
\begin{align}
  \hat{\pi}_\lambda(f_i)&=-x_i,
\end{align}
for $i=1,2,\dots,n$;
\item[2)]
\begin{align}
\begin{aligned}
  \pi_\lambda(h)&= E_{\hat{x}} + \sigma_{\lambda+\rho}(h), \\
  \pi_\lambda(h_A)&=- {\textstyle \sum\limits_{i,j=1}^n}
	a_{ij}\hat{x}_j\partial_{\hat{x}_i}	+\sigma_{\lambda+\rho}(h_A),
\end{aligned}
\end{align}
\begin{align}
\begin{aligned}
  \hat{\pi}_\lambda(h)&= -E_x+ \sigma_{\lambda-\rho}(h), \\
  \hat{\pi}_\lambda(h_A)&= {\textstyle \sum\limits_{i,j=1}^n}
	a_{ij}x_i\partial_{x_j}	+ \sigma_{\lambda-\rho}(h_A)
\end{aligned}
\end{align}
for $A \in M_{n, n}(\C)$ satisfying $A^{\rm t}I_{p,q}+I_{p,q}A=0$;
\item[3)]
\begin{align}
  \pi_\lambda(g_i)=-{\textstyle {1\over 2}} \veps_i {\textstyle \sum\limits_{j=1}^n} \veps_j \hat{x}_j^2 \partial_{\hat{x}_i}+\hat{x}_iE_{\hat{x}}
   +\hat{x}_i \sigma_{\lambda+\rho}(h)+ {\textstyle \sum\limits_{j=1}^n} \hat{x}_j
   \sigma_{\lambda+\rho}(h_{\veps_i\veps_j E_{ij}-E_{ji}})
\end{align}
\begin{align}
  \hat{\pi}_\lambda(g_i)=-{\textstyle {1\over 2}}\veps_i x_i {\textstyle \sum\limits_{j=1}^n} \veps_j\partial_{x_j}^2 + \partial_{x_i}E_x-
  \partial_{x_i}\sigma_{\lambda-\rho}(h) - {\textstyle \sum\limits_{j=1}^n} \partial_{x_j}
  \sigma_{\lambda-\rho}(h_{\veps_i\veps_j E_{ij}-E_{ji}})
\end{align}
for $i=1,2,\dots,n$.
\end{enumerate}}

\proof{The proof is a direct consequence of a straightforward but tedious
verification of all commutation relations. Another possibility is the
application of general formula for the representation action of $\mfrak{g}$ given in,
e.g.\ \cite{Krizka-Somberg2015}.}

Now, we shall fix a realization of the twisted complex (semi)simple spinor representation
$(\sigma_\lambda, \mS_\lambda)$  of $\mfrak{so}(p,q,\R)$.
Since the simple part $\mfrak{l}^{\rm s}$ of the complex Levi subalgebra $\mfrak{l}$ is isomorphic to
$\mfrak{so}(n,\C)$, we realize the (semi)simple
spinor module of $\mfrak{so}(n,\C)$ as the representation of $\mfrak{l}^{\rm s}$ on the exterior algebra of a Lagrangian
subspace in $\mC^n$. Let us denote by $\mS_{\pm}^n$ the irreducible half-spinor representations
for $\mfrak{l}^{\rm s} \simeq \mfrak{so}(n,\C)$ with $n$ even, and by
$\mS^{n}$ the spinor representation for $\mfrak{l}^{\rm s} \simeq \mfrak{so}(n,\C)$
with $n$ odd.
The generators of $\mfrak{l}^{\rm s}$ act on the spinor module
by the Clifford multiplication
\begin{align}
\begin{aligned} \label{eq:spinor-repr}
\mfrak{so}(n,\mC) \rarr \mcal{C}\ell_{p,q},
  \quad
  \sigma(h_{\veps_i\veps_jE_{ij}-E_{ji}})=
	-{\textstyle {1 \over 2}} \veps_i e_i e_j - {\textstyle {1\over 2}} \delta_{ij}
\end{aligned}
\end{align}
for all $i,j=1,2,\dots ,n$.
We used the convention $\mS=\mS_+^n \oplus \mS_-^n$ for $n$ even and
$\mS=\mS^n$ for $n$ odd, and denoted by $\mcal{C}\ell_{p,q}$ the
complex Clifford algebra for the symmetric bilinear form given
by $\langle v,w \rangle_{p,q}=v^{\rm t} I_{p,q} w$. The complex Clifford algebra $\mcal{C}\ell_{p,q}$
is an associative unital $\mC$-algebra given by quotient of the tensor
algebra $T(\C^n)$ by a two-sided ideal
$I \subset T(\C^n)$, generated by
\begin{align}
v\cdot w+w\cdot v=-2 \langle v,w \rangle_{p,q}1
\end{align}
for all $v,w \in\C^n$. In the canonical basis
$\{e_1,e_2,\dots,e_{p+q}\}$ of $\R^{p,q}$, $\R^{p,q}\otimes_\mR\mC\subset \mcal{C}\ell_{p,q}$,
we have $e_i\cdot e_i = -\veps_i 1$ for $i=1,2,\dots,p+q$.

The representation of $\mfrak{l}^{\rm s}$ extends to a representation
of $\mfrak{p}$ by the trivial action of the center $\mfrak{z}(\mfrak{l})$
of $\mfrak{l}$ and by the trivial action of the nilradical $\mfrak{u}$ of
$\mfrak{p}$. We retain the
notation $\sigma \colon \mfrak{p} \rarr \mfrak{gl}(\mathbb{S})$
for the extended action of the parabolic subalgebra
$\mfrak{p}$ of $\mfrak{g}$. In what follows, we are interested in the twisted
$\mfrak{p}$-module
$\sigma_\lambda \colon \mfrak{p} \rarr \mfrak{gl}(\mathbb{S}_\lambda)$
with a twist $\lambda \in \Hom_P(\mfrak{p},\C)$.
\medskip

The following $\mfrak{l}^{\rm s}$-invariant differential operators
\begin{align}
 D={\textstyle \sum\limits_{j=1}^n} e_j\partial_{x_j}, \qquad
 E={\textstyle \sum\limits_{j=1}^n} x_j\partial_{x_j}, \qquad
 X={\textstyle \sum\limits_{j=1}^n} \veps_je_j{x_j}
\end{align}
are the generators of the orthosymplectic Lie superalgebra
$\mfrak{osp}(1,2,\C)$. The consequences of our conventions for the complex Clifford algebra $\mcal{C}\ell_{p,q}$ include
\begin{align}
D^2= -{\textstyle \sum\limits_{j=1}^n} \veps_j \partial_{x_j}^2,\qquad
X^2= -{\textstyle \sum\limits_{j=1}^n} \veps_jx_j^2.
\end{align}
for $n\in \N$.

\medskip

\theorem{\label{thm:operator realization}
Let $\lambda \in \Hom_P(\mfrak{p},\C)$. Then the embedding of $\mfrak{g}$ into
$\eus{A}^\mfrak{g}_{\widebar{\mfrak{u}}^*}\! \otimes_\C \End \mathbb{S}_{\lambda-\rho}$
is given by
\begin{align}\label{ssw}
  \hat{\pi}_\lambda(f_i)=-x_i,
\end{align}
for $i=1,2,\dots,n$,
\begin{align}\label{spn}
\begin{aligned}
 \hat{\pi}_\lambda(h)&= -E_x + \lambda- {\textstyle {n \over 2}}, \\
   \hat{\pi}_\lambda(h_A)&= {\textstyle \sum\limits_{i,j=1}^n}
	a_{ij} x_i\partial_{x_j} + \sigma(h_A)
\end{aligned}
\end{align}
for $A \in M_{n, n}(\C)$ satisfying $A^{\rm t}I_{p,q}+I_{p,q}A=0$,
\begin{align}\label{soloper}
 \hat{\pi}_\lambda(g_i)={\textstyle {1\over 2}}\veps_i x_i D^2 + \partial_{x_i}(E_x- \lambda +{\textstyle {n\over 2}+ {1\over 2}})
   + {\textstyle {1\over 2}} \veps_i e_i D
\end{align}
for $i=1,2,\dots,n$.}

\proof{The proof is a straightforward combination of Theorem \ref{reproper} and the spinor
representation \eqref{eq:spinor-repr} twisted by character $\lambda -\rho \in \Hom_P(\mfrak{p},\C)$.
A direct derivation of this action for the twisted spinor representations is given in
\cite{KOSS}.}


\section{Generalized Verma modules and singular vectors}
\label{sec:singular vectors}

In what follows the generators $x_1,x_2,\dots,x_n$ of the graded
commutative $\C$-algebra $\C[\widebar{\mfrak{u}}^*]$ have the grading
$\deg(x_i)=1$ for $i=1,2,\dots,n$. As there is a canonical isomorphism of left $\eus{A}^\mfrak{g}_{\widebar{\mfrak{u}}^*}$-modules
\begin{align}
  \C[\widebar{\mfrak{u}}^*] \riso \eus{A}^\mfrak{g}_{\widebar{\mfrak{u}}^*}/\mcal{F}(I_e),
\end{align}
we obtain the isomorphism
\begin{align}
 \tau \circ \Phi_\lambda \colon  M^\mfrak{g}_\mfrak{p}(\mathbb{S}_{\lambda-\rho}) \riso
	\C[\widebar{\mfrak{u}}^*] \otimes_\C \mathbb{S}_{\lambda-\rho},
\end{align}
where the action of $\mfrak{g}$ on
$\C[\widebar{\mfrak{u}}^*] \otimes_\C \mathbb{S}_{\lambda-\rho}$
is given by Theorem \ref{thm:operator realization}.
Let us note that $\C[\widebar{\mfrak{u}}^*] \otimes_\C \mathbb{S}_{\lambda-\rho}$ and also
$\Sol(\mfrak{g},\mfrak{p}; \C[\widebar{\mfrak{u}}^*] \otimes_\C
\mathbb{S}_{\lambda-\rho})^\mcal{F}
\subset\C[\widebar{\mfrak{u}}^*] \otimes_\C \mathbb{S}_{\lambda-\rho}$ are semisimple
$\mfrak{l}$-modules.

It is well-known that the Fischer decomposition (cf.\ Appendix \ref{app:Fischer decompostion})
for the spinor-valued polynomials yields
an $\mfrak{l}$-module isomorphism
\begin{align}
  \varphi \colon \C[\widebar{\mfrak{u}}^*] \otimes_\C \mathbb{S}_{\lambda-\rho}
	\riso  \bigoplus_{a,b\in \N_0}  X^b M_a,
	\label{eq:decomposition C}
\end{align}
where $M_a = M_a^+ \oplus M_a^-$ for even $n$ and $M_a$ for odd $n$, respectively,
is the subspace of $\ker D$ of $a$-homogeneous $\mathbb{S}_{\lambda-\rho}$-valued
polynomials in the variables $(x_1,\dots,x_n)$ and
$X\in \eus{A}^\mfrak{g}_{\widebar{\mfrak{u}}^*}\!\otimes_\C \End \mathbb{S}_{\lambda-\rho}$ is
$\mfrak{l}^{\rm s}$-invariant.

\medskip

\lemma{\label{DXpowersrelations}
Let $k,m \in \N_0$. Then we have for all $v_m\in M_m$
\begin{enumerate}
\item[1)]
\begin{align}
  DX^kv_m = -k\, X^{k-1}v_m, \label{DXpowerseven}
\end{align}
if $k$ is even, and
\item[2)]
\begin{align}
  DX^kv_m= -(2m+n+k-1)\, X^{k-1}v_m, \label{DXpowersodd}
\end{align}
if $k$ is odd.
\end{enumerate}}

\proof{Using \eqref{comrel}, we have
\begin{align*}
DX^kv_m &= {\textstyle \sum\limits_{j=0}^{k-1}} (-1)^{j+1} X^j (2E+n) X^{k-j-1}v_m
= {\textstyle \sum\limits_{j=0}^{k-1}} (-1)^{j+1}(2(m+k-j-1)+n) X^{k-1}v_m
\end{align*}
for all $k,m\in\mN_0$ and $v_m \in M_m$. The specialization
to $k$ even and $k$ odd, respectively, implies the result.}

\lemma{Let $k \in \N_0$. Then we have
for all $j=1,2,\dots ,n$
\begin{enumerate}
\item[1)]
\begin{align}
[\partial_{x_j},X^k] = -k \veps_j x_j X^{k-2}, \label{eq:derivativexeven}
\end{align}
if $k$ is even, and
\item[2)]
\begin{align}
[\partial_{x_j},X^k] = \veps_j e_j X^{k-1}-(k-1)\veps_j x_j X^{k-2}, \label{eq:derivativexodd}
\end{align}
if $k$ is odd.
\end{enumerate}}

\proof{A direct computation gives $[\partial_{x_j},X]=\veps_j e_j$
and $[\partial_{x_j},X^2]=-2\veps_j x_j$ for all $j=1,2,\dots ,n$. Then for $k$ even, we have
\begin{align*}
[\partial_{x_j},X^k]={\textstyle \sum\limits_{r=0}^{{k\over 2}-1}} X^{2r}[\partial_{x_j},X^2] X^{k-2r-2} = -2 {\textstyle \sum\limits_{r=0}^{{k\over 2}-1}} \veps_j x_j X^{k-2} = -k\veps_j x_j X^{k-2}.
\end{align*}
For $k$ odd, we may write
\begin{align*}
[\partial_{x_j},X^k]=[\partial_{x_j},X]X^{k-1}+ X[\partial_{x_j},X^{k-1}]= \veps_j e_j X^{k-1}-(k-1)\veps_j x_jX^{k-2},
\end{align*}
where we used $[\partial_{x_j},X^{k-1}] = -(k-1) \veps_j x_j X^{k-3}$. The proof is complete.}

\medskip

\lemma{Let us introduce the differential operators $P_1,P_2,P_3 \in \eus{A}^\mfrak{g}_{\widebar{\mfrak{u}}^*}\!\otimes_\C \End \mathbb{S}_{\lambda-\rho}$ by
\begin{align}
P_1&={\textstyle \sum\limits_{j=1}^n} e_j\hat{\pi}_\lambda(g_j), \label{eq:sol operator 1}\\
P_2&={\textstyle \sum\limits_{j=1}^n} x_j\hat{\pi}_\lambda(g_j), \label{eq:sol operator 2}\\
P_3&={\textstyle \sum\limits_{j=1}^n} \veps_j \partial_{x_j}\hat{\pi}_\lambda(g_j). \label{eq:sol operator 3}
\end{align}
Then these operators are given by the explicit formulas
\begin{align}
P_1 &=\big(E_x-\lambda+{\textstyle {3\over 2}}+{\textstyle {1\over 2}}XD\big)D, \label{eq:sol operator 1 form}\\
P_2 &=-{\textstyle {1 \over 2}}X^2D^2+\big(E_x-\lambda+ {\textstyle {n \over 2}}+{\textstyle {1 \over 2}}\big)E_x+{\textstyle {1\over 2}} XD,
\label{eq:sol operator 2 form}\\
P_3 &=\big(\lambda-{\textstyle {1 \over 2}}E_x- 2\big)D^2.
\label{eq:sol operator 3 form}
\end{align}}

\proof{The proof is a direct consequence of \eqref{soloper} and the commutation relations
$[E_x,X]=X$, $[E_x,D]=-D$.}

Now, we shall find the subspace $\Sol(\mfrak{g},\mfrak{p},\C[\widebar{\mfrak{u}}^*] \otimes_\C
\mathbb{S}_{\lambda-\rho})^\mcal{F}$. Since $\Sol(\mfrak{g},\mfrak{p},\C[\widebar{\mfrak{u}}^*] \otimes_\C
\mathbb{S}_{\lambda-\rho})^\mcal{F}$ is a semisimple $\mfrak{l}$-module, we can assume that a solution $R$ of the system \eqref{soloper} is contained in some $\mfrak{l}$-isotypical component.
As $R \in \Sol(\mfrak{g},\mfrak{p},\C[\widebar{\mfrak{u}}^*] \otimes_\C
\mathbb{S}_{\lambda-\rho})^\mcal{F}$, we have $R \in \ker P_1 \cap \ker P_2 \cap \ker P_3$ by the construction of $P_1$, $P_2$ and $P_3$.

Therefore, we shall examine the common kernel of the differential
operators $P_1$, $P_2$ and $P_3$, relying on the results of Appendix \ref{app:Fischer decompostion}. In particular, we classify all solutions of the system of partial differential equations
given by $P_1, P_2, P_3$. After that we verify that these solutions are in the
solution space of \eqref{soloper} as well.

From Appendix \ref{app:Fischer decompostion} and \eqref{eq:decomposition C}, we know that the $\mfrak{l}$-isotypical components of $\C[\widebar{\mfrak{u}}^*] \otimes_\C \mathbb{S}_{\lambda-\rho}$ are of the form $X^kM^\pm_m$ for even $n$ and $X^k M_m$ for odd $n$, for all $k,m \in \N_0$. Let us assume that $R \in X^kM_m$, so we have $R=X^kv_m$ for $v_m \in M_m$.

If $k \in \N_0$ is even, then as a consequence of Lemma \ref{DXpowersrelations} we have
\begin{align}
  P_1X^kv_m&=-k\big({\textstyle {k \over 2}-\lambda-{n \over 2}+{3 \over 2}}\big)X^{k-1}v_m, \label{eq:operator p1 even}\\
  P_2X^kv_m&= \big({\textstyle (m+k)\big(m+k-\lambda+{n \over 2}+{1\over 2}\big)-{1 \over 2}k(2m+n+k-1)}\big) X^kv_m,  \label{eq:operator p2 even}\\
  P_3X^kv_m&= k(2m+n+k-2)\big({\textstyle \lambda - {1\over 2}(m+k+2)}\!\big)X^{k-2}v_m  \label{eq:operator p3 even}
\end{align}
for all $m \in \N_0$. From \eqref{eq:operator p1 even} we obtain that $k=0$ or $\lambda={k \over 2} -{n \over 2}+{3 \over 2}$. First of all, let us assume $k=0$. After substitution $k=0$ into \eqref{eq:operator p2 even} and \eqref{eq:operator p3 even}, we get either $m=0$, $\lambda \in \C$ or $m \neq 0$, $\lambda=m+{n\over 2}+{1 \over 2}$. Now, if we substitute $\lambda={k \over 2} - {n \over 2} +{3\over 2}$ into \eqref{eq:operator p2 even} and \eqref{eq:operator p3 even}, we obtain
\begin{align}
  P_2X^kv_m&= {\textstyle {1 \over 2}}(2m+k)(m+n-1) X^kv_m,  \\
  P_3X^kv_m&= -{\textstyle {1 \over 2}}k(2m+n+k-2)(m+n-1)X^{k-2}v_m,
\end{align}
which implies $k=0$, $m=0$, since we have $n \geq 3$.

Now, if $k \in \N_0$ is odd, then we have
\begin{align}
  P_1X^kv_m&=-(2m+n+k-1)\big({\textstyle {k \over 2}}-\lambda+m + 1\big)X^{k-1}v_m, \label{eq:operator p1 odd}\\
  P_2X^kv_m&= \big({\textstyle (m+k)\big(m+k-\lambda+{n \over 2}+{1\over 2}\big)-{1 \over 2}k(2m+n+k-1)}\big) X^kv_m,  \label{eq:operator p2 odd}\\
  P_3X^kv_m&= (k-1)(2m+n+k-1)\big({\textstyle \lambda - {1\over 2}(m+k+2)}\!\big)X^{k-2}v_m  \label{eq:operator p3 odd}
\end{align}
for all $m \in \N_0$. Therefore, from \eqref{eq:operator p1 odd} we obtain that $\lambda= m + {k \over 2} + 1$. Hence, after substitution $\lambda=m + {k \over 2} + 1$ into \eqref{eq:operator p2 odd} and \eqref{eq:operator p3 odd} we obtain
\begin{align}
  P_2X^kv_m&= {\textstyle {1\over 2}}m(n-k-1)X^kv_m,  \\
  P_3X^kv_m&= {\textstyle {1\over 2}}m(k-1)(2m+n+k-1)X^{k-2}v_m,
\end{align}
which implies $m=0$, since we have $n \geq 3$.

Therefore, there are three mutually exclusive cases giving potential
solutions of \eqref{soloper}:
\begin{enumerate}
  \item[1)] $m=0$, $k=0$, $\lambda \in \C$ and the $\mfrak{l}$-module is $M_0$;
  \item[2)] $m=0$, $k\in \N$ odd, $\lambda = {k \over 2}+ 1$
	and the $\mfrak{l}$-module is $X^kM_0$;
   \item[3)] $m \neq 0$, $k=0$, $\lambda=m+{\textstyle {n \over 2}+{1\over 2}}$ and the $\mfrak{l}$-module is $M_m$.
\end{enumerate}
We shall work out each case separately.
\bigskip

\noindent {\bf Case 1.} Let us assume $m=0$, $k=0$ and $\lambda \in \C$. Then we have $R=v_0$, where $v_0 \in M_0$. Since $v_0 \in \ker D$ and
$\partial_{x_i}v_0=0$ for $i=1,2,\dots,n$, we obtain
$\hat{\pi}_\lambda(g_i)v_0=0$ for $i=1,2,\dots,n$. Therefore, we have
\begin{align}
M_0 \subset \Sol(\mfrak{g},\mfrak{p},\C[\widebar{\mfrak{u}}^*] \otimes_\C
\mathbb{S}_{\lambda-\rho})^\mcal{F}
\end{align}
for $\lambda \in \C$.
\medskip

\noindent {\bf Case 2.} Let us assume $m=0$, $k \in \N$ odd and $\lambda = {k \over 2} + 1$. Then we have $R=X^kv_0$, where $v_0 \in M_0$. By Lemma \ref{DXpowersrelations} we have
\begin{align}
\begin{aligned}
DX^kv_0 &=-(n+k-1)X^{k-1}v_0, \\
D^2X^kv_0 &=(k-1)(n+k-1)X^{k-2}v_0
\end{aligned}
\end{align}
for all $v_0\in M_0$. Consequently, we get
\begin{align*}
\hat{\pi}_\lambda(g_i)X^kv_0&=
 \big({\textstyle {1\over 2}}\veps_i x_i D^2 + \partial_{x_i}(E_x- \lambda +{\textstyle {n\over 2}+{1\over 2}}) + {\textstyle {1\over 2}} \veps_i e_i D\big)X^kv_0 \\
& = {\textstyle {1 \over 2}}(k-1)(n+k-1)\veps_ix_iX^{k-2}v_0 \\
 & \quad + {\textstyle {1 \over 2}}(n+k-1)(X^k\partial_{x_i}+\veps_ie_iX^{k-1}-(k-1)\veps_ix_iX^{k-2})v_0 \\
& \quad - {\textstyle {1 \over 2}}(n+k-1)\veps_ie_iX^{k-1}v_0 \\
&= 0
\end{align*}
for all $v_0\in M_0$ and $i=1,2,\dots,n$.
Therefore, we have
\begin{align}
X^k M_0 \subset \Sol(\mfrak{g},\mfrak{p},\C[\widebar{\mfrak{u}}^*]
\otimes_\C \mathbb{S}_{\lambda-\rho})^\mcal{F}
\end{align}
for $\lambda = {k \over 2} + 1$ and $k$ odd natural number.
\medskip

\noindent {\bf Case 3.} Let us assume $m \neq 0$, $k=0$ and $\lambda = m + {n \over 2}+{1\over 2}$. Then we have $R=v_m$, where $v_m \in M_m$. Since $v_m \in \ker D$ and $E_xv_m=mv_m$, we obtain that $\hat{\pi}_\lambda(g_i)v_m=0$ for $i=1,2,\dots,n$. Therefore, we have
\begin{align}
M_m \subset \Sol(\mfrak{g},\mfrak{p},\C[\widebar{\mfrak{u}}^*]
\otimes_\C \mathbb{S}_{\lambda-\rho})^\mcal{F}
\end{align}
for $\lambda = m+{n \over 2}+{1\over 2}$ and $m\in \N$.
\medskip

\theorem{\label{thm:singular vector n odd}Let us assume $n\geq 3$ is odd. Then we have
\begin{align*}
  \tau \circ \Phi_{\lambda+\rho} \colon M^\mfrak{g}_\mfrak{p}(\mathbb{S}_\lambda)^\mfrak{u}
	\riso
  \begin{cases}
    M_0, & \text{if $\lambda-{1\over 2} \notin \N$, $\lambda + {n \over 2} -{1\over 2}\notin \N$}, \\
    M_0 \oplus M_{\lambda-{1 \over 2}}, & \text{if $\lambda - {1\over 2} \in \N$}, \\
    M_0 \oplus X^{2\lambda+n-2}M_0, & \text{if $\lambda+ {n\over 2}-{1\over 2} \in \N$}.
  \end{cases}
\end{align*}}

\proof{The proof follows from the discussion in Case 1 up to Case 3.}

\theorem{\label{thm:singular vector n even}Let us assume $n\geq 3$ is even. Then we have
\begin{align*}
  \tau \circ \Phi_{\lambda+\rho} \colon M^\mfrak{g}_\mfrak{p}(\mathbb{S}_\lambda)^\mfrak{u}
	\riso
  \begin{cases}
    M_0, & \text{if $\lambda + {n \over 2}-{1\over 2} \notin \N$}, \\
    M_0 \oplus X^{2\lambda+n-2}M_0, & \text{if $\lambda + {n \over 2}-{1\over 2} \in \N$, $\lambda -{1\over 2}\notin \N$}, \\
    M_0 \oplus M_{\lambda-{1 \over 2}} \oplus X^{2\lambda+n-2}M_0, & \text{if $\lambda-{1 \over 2} \in \N$}.
  \end{cases}
\end{align*}}

\proof{The proof follows from the discussion in Case 1 up to Case 3.}


\section{Equivariant differential operators on spinors
in conformal geometry}
\label{sec:invariant operators}

Given a complex (semi)simple finite-dimensional $\widetilde{P}$-module $(\sigma,\mathbb{V})$,
we consider the induced representation of $\smash{\widetilde{G}}$
on the space $\smash{\Ind_{\widetilde{P}}^{\widetilde{G}}}(\mV)$ of smooth sections
of the homogeneous vector bundle
$\mcal{V}=\smash{\widetilde{G}} \times_{\widetilde{P}} \mV \rarr \smash{\widetilde{G}}/\smash{\widetilde{P}}$,
\begin{align}
{\lC}^\infty(\widetilde{G}/\widetilde{P},\mcal{V})
\simeq{\lC}^\infty(\widetilde{G},\mV)^{\widetilde{P}}=\{f\in{\lC}^\infty(\widetilde{G},\mV);\,
f(gp)=\sigma(p^{-1})f(g),\ \text{for all}\ g\in \widetilde{G},\, p\in \widetilde{P}\}.
\end{align}
We denote by $J^k_e(\widetilde{G},\mV)^{\widetilde{P}}$ the space of $k$-jets in $e\in \widetilde{G}$
of $\widetilde{P}$-equivariant smooth mappings for $k\in \N_0$, and
by $\smash{J^\infty_e(\widetilde{G},\mV)^{\widetilde{P}}}$ its projective limit
\begin{align}
J^\infty_e(\widetilde{G},\mathbb{V})^{\widetilde{P}}
= {\textstyle \lim\limits_{\longrightarrow k}} J^k_e(\widetilde{G},\mV)^{\widetilde{P}}.
\end{align}
Then there is a non-degenerate $(\mfrak{g},\smash{\widetilde{P}})$-invariant
pairing between $\smash{J^\infty_e(\widetilde{G},\mV)^{\smash{\widetilde{P}}}}$  and
$M^\gog_\gop(\mV^*),$ which identifies the generalized Verma module
$M^\gog_\gop(\mV^*)$ with the vector space of all $\C$-linear mappings
$\smash{J^\infty_e(\widetilde{G},\mV)^{\widetilde{P}}}\rarr \mC$ that factor through
$\smash{J^k_e(\widetilde{G},\mathbb{V})^{\widetilde{P}}}$
for some $k \in \N_0$. Here we denoted by $\mV^*$ the dual
$\mathfrak{g}$-module equipped with the dual action of $\mathfrak{g}$.

There is a classical consequence of the last statement
explaining the relationship between the geometrical problem of finding
${\widetilde G}$-equivariant differential operators between induced representations
and the algebraic problem of finding homomorphisms between generalized
Verma modules, cf.\ \cite{Collingwood-Shelton1990}, \cite{Kostant1975}.
Let $\mathbb{V}$ and $\mathbb{W}$ be complex (semi)simple finite-dimensional
$\smash{\widetilde{P}}$-modules.
Then the vector space of $\smash{\widetilde{G}}$-equivariant differential operators
from $\smash{\Ind_{\widetilde{P}}^{\widetilde{G}}}(\mathbb{V})$
to $\smash{\Ind_{\widetilde{P}}^{\widetilde{G}}}(\mathbb{W})$ is isomorphic to
the vector space of $(\gog,\smash{\widetilde{P}})$-homomorphisms of generalized Verma modules
from $M^\gog_\gop(\mathbb{W}^*)$ to $M^\gog_\gop(\mathbb{V}^*)$.
\medskip

In the following theorem we retain the notation of Section \ref{sec:Verma modules};
for the group $\Spin(p,q,\R)$ we denote by $\mS$ the
spinor representation for $n=p+q$ odd, and the direct sum of the half-spinor
representations $\mS_\pm$ for $n=p+q$ even. If a character $\lambda \in
\Hom_{\widetilde P}(\mfrak{p},\C)$ defines a group character $e^\lambda$ of $\smash{\widetilde{P}}$, then the representation $\mS$ may be twisted by $e^\lambda$ to the representations $\mS_\lambda$ of the parabolic subgroup $\smash{\widetilde{P}}=\smash{\widetilde{M}}AU$,
$\smash{\widetilde{L}}=\smash{\widetilde{M}}A$, with $\smash{\widetilde{M}} \simeq \Spin(p,q,\mR)$ and
$A$ acting in the one-dimensional representation $\C_\lambda$ ($U$ acts trivially).
\medskip

\theorem{\label{thm:edo}
Let $\smash{\widetilde{G}}=\Spin_o(p+1,q+1,\R)$ be the identity component
of the spin group $\Spin(p+1,q+1,\R)$, $n=p+q\geq 3$ and $\lambda\in \Hom_{\widetilde{P}}(\mfrak{p},\C)$.
Furthermore, let $\smash{\widetilde{P}} \simeq (\GL(1,\R)_+\times \Spin(p,q,\R))\ltimes\R^{p,q}$ be the
maximal (conformal) parabolic subgroup of $\smash{\widetilde{G}}$ with the unipotent
radical in the Langlands-Iwasawa decomposition for $\smash{\widetilde{P}}$ isomorphic to $\R^{p,q}$.
For $\mathbb{V}= \mS_\lambda$ we have $\mathbb{V}^* \simeq \mS^*_{-\lambda}$. Then the singular vectors in Theorem \ref{thm:singular vector n odd} and Theorem \ref{thm:singular vector n even} correspond, in the non-compact picture of the induced representations, to $\smash{\widetilde{G}}$-equivariant differential operators as follows:
\begin{enumerate}
\item[1)] Let $\lambda=-{1\over 2}(n-2-a)\widetilde{\omega}$
and $\mu=-{1 \over 2}(n-2+a)\widetilde{\omega}$ for $a \in \N$
odd.
Then there are $\smash{\widetilde{G}}$-equivariant differential operators
\begin{align}
D_a \colon {\lC}^\infty(\widebar{\mfrak{u}}(\R),\mS^*_{-\lambda})
\rarr  {\lC}^\infty(\widebar{\mfrak{u}}(\R),\mS^*_{-\mu})
\end{align}
 of order $a \in \N$, $\widebar{\mfrak{u}}(\R)\simeq \mR^{p,q}$.
The infinitesimal intertwining property of $D_a$ is
\begin{align}
D_a\pi^*_{-{1 \over 2}(a+2)}(X)=\pi^*_{{1 \over 2}(a-2)}(X)D_a
\end{align}
for all $X \in \mfrak{g}$. We call these operators conformal powers of the Dirac operator.
\item[2)] Let $\lambda=\big(a+{1\over 2}\big)\widetilde{\omega}$ for $a\in \N$. Then
there are $\smash{\widetilde{G}}$-equivariant differential operators
\begin{align}
T_a \colon {\lC}^\infty(\widebar{\mfrak{u}}(\R),\mS^*_{-\lambda})
\rarr  {\lC}^\infty(\widebar{\mfrak{u}}(\R),(M_a)^*)
\end{align}
 of order $a \in \N$. These operators are called conformal twistor operators on spinors.
\end{enumerate}}

The remaining collection of singular vectors are \smash{$\widetilde{G}$}-equivariant differential
operator given by a multiple of the identity map.

\begin{appendices}

\section{The Fischer decomposition for $\mfrak{so}(n,\C)$}
\label{app:Fischer decompostion}

Here we recall a well-known result, cf.\ \cite{Delanghe-Sommen-Soucek1992}, describing the action
of the simple part of the Levi subalgebra
$\mfrak{l}^{\rm s}\simeq \mfrak{so}(n,\C)$-module structure on the
conformal generalized Verma module $M^\gog_\gop(\mS_{\lambda-\rho})$
for arbitrary $\lambda\in\C$. 

The decomposition is given by the spin Howe duality for the pair
$\mathfrak{osp}(1,2,\C) \oplus \mathfrak{so}(n,\C)$, where
$\mathfrak{so}(n,\C)$ acts on
$\C[\widebar{\mfrak{u}}^*] \otimes_\C \mathbb{S}_{\lambda-\rho}$
by \eqref{spn} and $\mathfrak{osp}(1,2,\C)$
acts by
\begin{align}
 D=\sum\limits_{j=1}^n \, e_j\partial_{x_j}, \qquad
 E=\sum_{j=1}^{n} \, x_{j}\partial_{x_j}, \qquad
 X=\sum\limits_{j=1}^n \, \veps_je_jx_j .
\end{align}
This decomposition is of the form
\begin{align}
\C[\widebar{\mfrak{u}}^*]\otimes_\C \mS_{\lambda-\rho}
\simeq \smash{\bigoplus_{a,b \in \N_0}} X^b M_a,
\quad M_a=\big(\C[\widebar{\mfrak{u}}^*]_a\otimes_\C \mS_{\lambda-\rho}\big)
\cap\ker D ,
\end{align}
\begin{align*}
\xymatrix@=11pt{P_0 \otimes \mathbb{S} \ar@{=}[d] &  P_1 \otimes \mathbb{S} \ar@{=}[d]&
P_2 \otimes \mathbb{S} \ar@{=}[d] & P_3 \otimes \mathbb{S} \ar@{=}[d] &
P_4 \otimes \mathbb{S} \ar@{=}[d]& P_5 \otimes \mathbb{S} \ar@{=}[d] & \ldots \\
M_0 \ar[r] & X M_0 \ar @{} [d] |{\oplus} \ar[r] & X^2 M_0 \ar @{} [d] |{\oplus} \ar[r] & X^3 M_0 \ar @{} [d] |{\oplus}
 \ar[r] & X^4 M_0 \ar @{} [d] |{\oplus}\ar[r] & X^5 M_0 \ar @{} [d] |{\oplus} & \ldots \\
& M_1 \ar[r] & X M_1 \ar @{} [d] |{\oplus}\ar[r] & X^2 M_1 \ar @{} [d] |{\oplus}
 \ar[r] & X_s^3 M_1 \ar @{} [d] |{\oplus}\ar[r] & X^4 M_1 \ar @{} [d] |{\oplus} & \ldots \\
&& M_2 \ar[r] & X M_2 \ar @{} [d] |{\oplus}
 \ar[r] & X^2 M_2 \ar @{} [d] |{\oplus}\ar[r] & X^3 M_2 \ar @{} [d] |{\oplus} & \ldots \\
&&& M_3 \ar[r] & X M_3 \ar @{} [d] |{\oplus}\ar[r] & X^2 M_3  \ar @{} [d] |{\oplus} & \ldots \\
&&&& M_4 \ar[r] & X M_4 \ar @{} [d] |{\oplus} & \ldots  \\
&&&&& M_5 & \ldots}
\end{align*}
In the scheme above we used the notation $P_a=\C[\widebar{\mfrak{u}}^*]_a$
for the $a$-homogeneous polynomials and $M_a=M_a^+ \oplus M_a^-$ for even $n$ and $M_a$ for odd $n$. The operators $D$ and $X$ act in
the previous picture horizontally, $E$ preserves each simple orthogonal module
in the decomposition and the $\mathfrak{osp}(1,2,\C)$-commutation relations
of $\mathfrak{so}(n,\C)$-invariant operators are
\begin{align}\label{comrel}
[E,D]=-D, \qquad \{ D,X\}= -2E-n, \qquad [E,X]= X.
\end{align}

This implies the Fischer decomposition for any real form $\mfrak{so}(p,q,\R)$ of $\mfrak{so}(n,\C)$ on the
space of spinor valued polynomials on the real subspace $\mR^{p,q}$ of $\C^n$,
$p+q=n$. In more detail, since $M_a^+$ and $M_a^-$ for even $n$ and $M_a$ for odd $n$ are irreducible representations of $\mfrak{so}(p,q,\R)$, and the operators $E,X,D$ are $\mfrak{so}(p,q,\R)$-invariant, we obtain that $X^bM_a$ are representations of $\mfrak{so}(p,q,\R)$ for all $a,b \in \N_0$.
The action of $\mfrak{so}(p,q,\R)$ on $\R[\widebar{\mfrak{u}}(\R)^*] \otimes_\R \mathbb{S}_{\lambda-\rho}$ is given by
\begin{align}
  \veps_i\veps_jE_{ij}-E_{ji} \mapsto  \veps_i\veps_j x_i\partial_{x_j}-x_j \partial_{x_i} - {\textstyle {1 \over 2}}\veps_i e_i e_j
\end{align}
for all $i,j=1,2,\dots,n$ satisfying $i \neq j$. Because the relations \eqref{comrel}
are independent of the signature, the Fischer decomposition does not depend on the 
signature of the real form as well.

\end{appendices}


\section*{Acknowledgement}

L.\,Křižka is supported by PRVOUK p47,
P.\,Somberg acknowledges the financial support from the grant GA\,CR P201/12/G028.



\begin{thebibliography}{10}

\bibitem{Baum-Friedrich-Grunewald-Kath1991}
Helga Baum, Thomas Friedrich, Ralf Grunewald, and Ines Kath, \emph{{Twistors
  and Killing spinors on Riemannian manifolds}}, Teubner-Texte zur Mathematik,
  vol. 124, Teubner, Stuttgart, 1991.

\bibitem{Boe-Collingwood1985}
Brian~D. Boe and David~H. Collingwood, \emph{{A comparison theory for the
  structure of induced representations. II}}, Math. Z. \textbf{190} (1985),
  no.~1, 1--11.

\bibitem{Collingwood-Shelton1990}
David~H. Collingwood and Brad Shelton, \emph{{A duality theorem for extensions
  of induced highest weight modules}}, Pacific J. Math. \textbf{146} (1990),
  no.~2, 227--237.

\bibitem{Curry-Gover2014}
Sean~N. Curry and Rod~A. Gover, \emph{{An introduction to conformal geometry
  and tractor calculus, with a view to applications in general relativity}},
  {\tt arXiv:1412.7559} (2014).

\bibitem{Delanghe-Sommen-Soucek1992}
Richard Delanghe, Franciscus Sommen, and Vladimír Souček, \emph{{Clifford
  algebra and spinor-valued functions. A function theory for the Dirac
  operator}}, Mathematics and its Applications, vol.~53, Kluwer Academic
  Publishers Group, Dordrecht, 1992.

\bibitem{Eastwood-Slovak1997}
Michael Eastwood and Jan Slovák, \emph{{Semiholonomic Verma modules}}, J.
  Algebra \textbf{197} (1997), no.~2, 424--448.

\bibitem{Fischmann2014}
Matthias Fischmann, \emph{{On conformal powers of the Dirac operator on spin
  manifolds}}, Arch. Math. \textbf{50} (2014), no.~4, 237--254.

\bibitem{Friedrich2000}
Thomas Friedrich, \emph{{Dirac operators in Riemannian geometry}}, Graduate
  Studies in Mathematics, vol.~25, American Mathematical Society, Providence,
  2000.

\bibitem{Holland-Sparling2001}
Jonathan Holland and George Sparling, \emph{{Conformally invariant powers of
  the ambient Dirac operator}}, {\tt arXiv:math/0112033} (2001).

\bibitem{Krizka-Somberg2015}
Libor Křižka and Petr Somberg, \emph{{Algebraic analysis of scalar generalized
  Verma modules of Heisenberg parabolic type I.: $A_n$-series}}, {\tt
  arXiv:1502.07095} (2015), accepted for publication in Transformation Groups.

\bibitem{KOSS}
Toshiyuki Kobayashi, Bent \O{}rsted, Petr Somberg, and Vladimír Souček,
  \emph{{Branching laws for Verma modules and applications in parabolic
  geometry. I}}, Adv. Math. \textbf{285} (2015), 1796--1852.

\bibitem{Kostant1975}
Bertram Kostant, \emph{{Verma modules and the existence of quasi-invariant
  differential operators}}, {Non-Commutative Harmonic Analysis}, Lecture Notes
  in Mathematics, vol. 466, Springer, Berlin, 1975, pp.~101--128.

\bibitem{Slovak1992}
Jan Slovák, \emph{{Natural operators on conformal manifolds}}, Research Lecture
  Notes, University of Vienna, 1992, {(Habilitation Thesis, Masaryk University,
  1993)}.

\end{thebibliography}

\providecommand{\bysame}{\leavevmode\hbox to3em{\hrulefill}\thinspace}
\providecommand{\MR}{\relax\ifhmode\unskip\space\fi MR }
\providecommand{\MRhref}[2]{%
  \href{http://www.ams.org/mathscinet-getitem?mr=#1}{#2}
}
\providecommand{\href}[2]{#2}

\end{document}